\documentclass[11pt]{article}

\usepackage[T2A]{fontenc}
\usepackage[cp1251]{inputenc}
\usepackage{amsthm}
\usepackage{amsfonts}
\usepackage{eufrak}
\usepackage{amssymb}
\usepackage{amsmath}
\usepackage{cite}


\begin{document}
\newtheoremstyle{mytheorem}
  {\topsep}   
  {\topsep}   
  {\itshape}  
  {}       
  {\bfseries} 
  {  }         
  {5pt plus 1pt minus 1pt} 
  { }          
\newtheoremstyle{myremark}
  {\topsep}   
  {\topsep}   
  {\upshape}  
  {}       
  {\bfseries} 
  {  }         
  {5pt plus 1pt minus 1pt} 
  { }          
\theoremstyle{mytheorem}
\newtheorem{theorem}{Theorem}[section]
 \newtheorem{theorema}{Theorem}
 \newtheorem*{a*}{The Kac--Bernstein theorem}
 \newtheorem*{b*}{Theorem A}
\newtheorem{proposition}{Proposition}[section]
\newtheorem{lemma}{Lemma}[section]
\newtheorem{corollary}{Corollary}[section]
\newtheorem{definition}{Definition}[section]
\theoremstyle{myremark}
\newtheorem{remark}{Remark}[section]

\centerline{\textbf{Solution of the Kac--Bernstein functional equation }}

\centerline{\textbf{  on   Abelian groups in the class of positive functions}}

\bigskip

\centerline{\bf G.M. Feldman}

\bigskip

\makebox[20mm]{ }\parbox{125mm}{ \small   The general form of the solutions of the   Kac--Bernstein functional equation
$$
f(x+y)g(x-y)=f(x)f(y)g(x)g(-y), \ x, y\in X,
$$
on an arbitrary Abelian  group $X$ in the class of positive functions is obtained.  We also study the solutions of this equation in the class of complex-valued functions that do not vanish and satisfy the Hermitian condition.}

\bigskip


\bigskip


\section{ Introduction}

Let $X$ be an Abelian group. We remaind that  the Kac--Bernstein functional equation on the group $X$ is an equation of the form
\begin{equation}\label{1}
f(x+y)g(x-y)=f(x)f(y)g(x)g(-y), \ x, y\in X.
\end{equation}
The Kac--Bernstein functional equation arises in the proof of the following characterization theorem of mathematical statistics on the real line.
\begin{a*}  Let  $\xi_1$ and $\xi_2$ be independent random variables.   If the sum $\xi_1+\xi_2$  and the difference $\xi_1-\xi_2$ are also independent, then   $\xi_1$ and $\xi_2$ are Gaussian random variables.
\end{a*}

This theorem was independently proved by M. Kac  \cite{Ka} and S.N. Bernstein    \cite{Be}. The Kac--Bernstein theorem   is one of the first, and  apparently, the most famous characterization theorem of mathematical statistics. Some generalizations of this theorem were studied in the case when independent random variables take values either in an arbitrary locally compact Abelian group or in specific locally compact Abelian groups
(\cite{BaESta,  Fe5, Fe15, TVP2011, TVP2016,  HeRa, Mi1, Mi3, Ru1, Ru2},
see also \cite[Chapter 3]{Fe5n}). In this case, either the densities of the distributions of independent random variables or their characteristic functions satisfy the  equation (\ref{1}).  In the latter case  the equation (\ref{1}) is considered on the character group of the original group. We also remark that the   equation (\ref{1}) was studied in \cite{FeM1} on locally compact Abelian groups for complex-valued continuous and measurable functions that satisfy the Hermitian condition.

In the article \cite{Ka1} Pl. Kannappan  considered the equation  (\ref{1}) on the real line assuming that the functions  $f(x)$ and $g(x)$ are positive. Moreover, no restrictions were imposed on the functions $f(x)$ and $g(x)$   such as continuity, measurability, positive definiteness.
Pl. Kannappan proved the following theorem.
\begin{b*} \label{K}
Let $f(x)$ and $g(x)$  be positive functions on the real line satisfying the Kac--Bernstein functional equation $(\ref{1})$. Then the functions $f(x)$ and $g(x)$ may be represented in the form
$$
f(x)=\exp\{P(x)+l(x)+r\}, \quad g(x)=\exp\{P(x)+m(x)-r\},
$$
where the function $P(x)$ satisfies the equation
\begin{equation}\label{2}
P(x+y)+P(x-y)=2[P(x)+P(y)],
\end{equation}
each of the functions $l(x)$ and $m(x)$ satisfies the Cauchy equation
\begin{equation}\label{3}
l(x+y)=l(x)+l(y),
\end{equation}
$r$ is a constant.
\end{b*}
The main result of this note is the proof of a theorem in which a complete description of the solutions of the equation (\ref{1}) in the class of positive functions on an arbitrary Abelian group $X$ is given.

\section{Solution of the Kac--Bernstein functional equation in the class of positive functions}

Let $X$ be an Abelian group, $n$ be a positive integer. Denote by $X^{(n)}$ the image of the group $X$ under the mapping $x\rightarrow nx$. Let
$$
X=\bigcup_\alpha\left(x_\alpha+X^{(2)}\right)
$$
be a decomposition of the group $X$ with respect to the subgroup $X^{(2)}$, where an index $\alpha$ runs through the factor-group $X/X^{(2)}$. Moreover, we assume that the element $x_0=0$ corresponds  to the coset $X^{(2)}$.

Let $T(x)$ be a function on the group $X$, and $h$ be an arbitrary element of the group   $X$. Denote by $\Delta_h$ the finite difference operator
$$
\Delta_hT(x)=T(x+h)-T(x).
$$
A function $T(x)$ is called a polynomial on the group  $X$ if for some nonnegative integer   $n$ the function $T(x)$ satisfies the equation
\begin{equation}\label{4}
\Delta_h^{n+1}T(x)=0, \quad x, h\in X.
\end{equation}
The minimal $n$ at which this equation holds is called the degree of the polynomial  $T(x)$. We remark that if a function $P(x)$ satisfies the equation $(\ref{2})$, then either  $P(x)\equiv 0$ for $x\in X$ or $P(x)$ is a polynomial of degree 2. If a function $l(x)$ satisfies the equation  $(\ref{3})$, then either  $l(x)\equiv 0$ for $x\in X$ or $l(x)$ is a polynomial of degree 1.

The main result of the article is the proof of the  following statement.
\begin{theorem}\label{th1}
Let $X$ be an Abelian group. Let $f(x)$ and $g(x)$ be positive functions on the group    $X$ satisfying the Kac--Bernstein functional equation   $(\ref{1})$. Then the functions $f(x)$ and $g(x)$ can be represented in the form
\begin{equation}\label{rev1}
f(x)=\exp\{P(x)+l(x)+r(x)\},  \quad
 g(x)=\exp\{P(x)+m(x)-r(x)\}, \quad x\in X,
\end{equation}
where the function $P(x)$ satisfies the equation  $(\ref{2})$,
each of the functions $l(x)$ and $m(x)$ satisfies the   equation  $(\ref{3})$,
the function $r(x)$ is  a constant    on each of the cosets   of the subgroup $X^{(2)}$ in  $X$.
\end{theorem}
To prove Theorem  \ref{th1} we need the following lemmas.
\begin{lemma}\label{l1} {\rm(\!\cite{Djo})}
Let $X$ be an Abelian group, $T(x)$  be a function on the group $X$  satisfying the equation   $(\ref{4})$. Then there are symmetric $k$-additive functions
$g_k(x_1, \dots, x_k)$ on the group $X^k$, $k=1, 2, \dots, n$, such that
$$
T(x)=\sum_{k=0}^ng_k^{*}(x), \quad x\in X,
$$
where $g^{*}_0(x)=T(0)$, $g^{*}_k(x)=g_k(x, \dots, x)$, $k=1, 2, \dots, n$.
\end{lemma}
\begin{corollary}\label{c1}
Let $X$ be an Abelian group, $T(x)$  be a polynomial on the group   $X$ of degree $\le 2$. Then $T(x)$ may be represented in the form
$$
T(x)=P(x)+l(x)+r, \quad x\in X,
$$
where the function $P(x)$ satisfies the equation    $(\ref{2})$,
  the function $l(x)$ satisfies the equation $(\ref{3})$,
$r$ is a constant.
\end{corollary}
\begin{lemma}\label{l2}
Let $X$ be an Abelian group, $T(x)$  be a function on the group $X$ satisfying the equation
\begin{equation}\label{5}
\Delta_{2k}\Delta_h^{2}T(x)=0, \quad k, h, x\in X.
\end{equation}
Then the function  $T(x)$   may be represented in the form
\begin{equation}\label{6}
 T(x)=P(x)+l(x)+r(x), \quad x\in X,
\end{equation}
where the function $P(x)$ satisfies the equation    $(\ref{2})$,
  the function $l(x)$ satisfies the equation $(\ref{3})$,
    the function $r(x)$ is  a constant    on each of the cosets of the subgroup $X^{(2)}$ in  $X$.
\end{lemma}
{\bf Proof}. Note that the equation (\ref{5}) on specific groups  $X$ for an even function $T(x)$   arose  in solving characterization problems of mathematical statistics on locally compact Abelian groups. Lemma \ref{l2} on an arbitrary locally compact abelian group $X$ for a continuous even function $T(x)$ was proved in  \cite[Lemma 6]{Fe20bb}. We will follow the scheme of the proof of Lemma 6 in \cite{Fe20bb}.

Obviously, without loss of generality, we can assume that $T(0)=0$. We fix an element $x_\alpha\in X$. It follows from the equation (\ref{5}) that the restriction of the function   $T(x_\alpha+x)$ to the subgroup $X^{(2)}$ is a polynomial of degree $\le 2$. Therefore, by Corollary \ref{c1} applied to the subgroup   $X^{(2)}$,  there is a representation
\begin{equation}\label{7}
 T(x_\alpha+x)=P_\alpha(x)+l_\alpha(x)+s_\alpha, \quad x\in X^{(2)},
\end{equation}
where the function $P_\alpha(x)$ satisfies the equation    $(\ref{2})$,
  the function $l_\alpha(x)$ satisfies the equation  $(\ref{3})$,
$s_\alpha$ is a constant. In particular, $s_0=0$. Moreover,   by Lemma \ref{l1}, $P_\alpha(x)=A_\alpha(x, x)$, where $A_\alpha(x_1, x_2)$, $x_1, x_2\in X^{(2)}$, is a symmetric biadditive   function. We extend the function $A_\alpha(x_1, x_2)$ from the subgroup $(X^{(2)})^2$ as  a symmetric biadditive   function    $\tilde A_\alpha(x_1, x_2)$ to the group $X^2$, setting $\tilde A_\alpha(x_1, x_2)={1\over 4}A_\alpha(2x_1, 2x_2)$, $x_1, x_2\in X$. Denote by $\tilde P_\alpha(x)$ the corresponding extension of the function $P_\alpha(x)$ from $X^{(2)}$ to $X$. We extend the function $l_\alpha(x)$ from the subgroup $X^{(2)}$ as an additive function $\tilde l_\alpha(x)$ to the group $X$, setting   $\tilde l_\alpha(x)={1\over 2}l_\alpha(2x)$, $x\in X$.   It follows from $(\ref{7})$ that
\begin{equation}\label{8}
 T(-x_\alpha-x)=T(x_\alpha+(-2x_\alpha-x)
 =P_\alpha(-2x_\alpha-x)+l_\alpha(-2x_\alpha-x)+s_\alpha$$$$=4\tilde P_\alpha(x_\alpha)+4\tilde A_\alpha(x_\alpha, x)+P_\alpha(x)-2\tilde l_\alpha(x_\alpha)- l_\alpha(x)+s_\alpha, \quad x\in X^{(2)}.
\end{equation}
Substituting $k=h$ in  $(\ref{5})$, we obtain
\begin{equation}\label{9}
 T(x+4h)-2T(x+3h)+2T(x+h)-T(x)=0, \quad x, h\in X.
\end{equation}
Substituting $x=0$, $h=x_\alpha+2x$  in (\ref{9}), we find
\begin{equation}\label{10}
 T(4x_\alpha+8x)-2T(3x_\alpha+6x)+2T(x_\alpha+2x)=0, \quad x\in X.
\end{equation}

 1. Suppose that $T(-x)=T(x)$, i.e. $T(x)$ is an even function. Then  (\ref{7}) and (\ref{8}) imply that
\begin{equation}\label{11}
 l_\alpha(x)=2\tilde A_\alpha(x_\alpha, x), \quad x\in X^{(2)}.
\end{equation}
Since $P_\alpha(x)=\tilde  P_\alpha(x_\alpha+x)-2\tilde A_\alpha(x_\alpha, x)-\tilde  P_\alpha(x_\alpha)$, we find from   (\ref{7}) and (\ref{11}) that
\begin{equation}\label{12}
T(x_\alpha+x)=\tilde P_\alpha(x_\alpha+x)+r_\alpha, \quad x\in X^{(2)},
\end{equation}
where $r_\alpha=-\tilde P_\alpha(x_\alpha)+s_\alpha$. We note that if $\alpha=0$,  then   $r_0=0$. Taking into account (\ref{12}), we obtain from  (\ref{10})
$$
\tilde P_0(4x_\alpha+8x)-2\tilde P_\alpha(3x_\alpha+6x)-2r_\alpha+2\tilde P_\alpha(x_\alpha+2x)+2r_\alpha=0, \quad x\in X.
$$
This implies that
\begin{equation}\label{13}
4(\tilde P_0(x)-\tilde P_\alpha(x))+4(\tilde A_0(x_\alpha, x)-\tilde A_\alpha(x_\alpha, x))+\tilde P_0(x_\alpha)-\tilde P_\alpha(x_\alpha)=0, \quad x\in X.
\end{equation}
It follows from (\ref{13}) that $\tilde P_0(x)=\tilde P_\alpha(x)$ for all $x\in X$ and any $\alpha$. Taking this into account, (\ref{12}) implies for the function  $T(x)$ the representation
\begin{equation}\label{n1}
T(x)=\tilde P_0(x)+r_\alpha, \quad x\in x_\alpha+X^{(2)}.
\end{equation}
Put
\begin{equation}\label{n6}
P(x)=\tilde P_0(x), \quad x\in X, \quad r(x)=r_\alpha, \quad x\in x_\alpha+X^{(2)}.
\end{equation}
The desired representation for the function $T(x)$ follows from  (\ref{n1}) and (\ref{n6}). Moreover, in the case when $T(x)$ is an even function, we proved that    $l(x)\equiv 0$ for $x\in X$ in the representation (\ref{6}).

2. Suppose that $T(-x)=-T(x)$, i.e. $T(x)$ is an odd function.   Then  (\ref{7}) and (\ref{8}) imply that
\begin{equation}\label{14}
P_\alpha(x)=0, \quad x\in X^{(2)},  \quad \tilde l_\alpha(x_\alpha)=s_\alpha.  \end{equation}
Taking into account (\ref{14}), it follows from (\ref{7})    that
\begin{equation}\label{15}
T(x_\alpha+x)=\tilde l_\alpha(x_\alpha+x), \quad x\in X^{(2)}.
\end{equation}
Hence, the equation (\ref{10}) takes the form
$$
\tilde l_0 (4x_\alpha+8x)-2\tilde l_\alpha(3x_\alpha+6x)+2\tilde l_\alpha(x_\alpha+2x)=0, \quad x\in X.
$$
This implies that $\tilde l_0(x)=\tilde l_\alpha(x)$ for $x\in X$. In view of this,   (\ref{15}) implies that
\begin{equation}\label{n2}
T(x)=\tilde l_0(x), \quad x\in X.
\end{equation}
Put $l(x)=\tilde l_0(x)$, $x\in X.$ The desired representation for the function $T(x)$ follows from  (\ref{n2}). Moreover, in the case when $T(x)$ is an odd function, we proved that $P(x)\equiv 0$, $r(x)\equiv 0$ for $x\in X$  in the representation (\ref{6}).

To complete the proof of the lemma we note that  an arbitrary function $T(x)$   is   the sum of an even and an odd function. $\blacksquare$
\medskip

{\bf Proof of Theorem \ref{th1}}. Put $T(x)=\log f(x)$ and $S(x)=\log g(x)$. It follows from (\ref{1}) that the functions $T(x)$ and $S(x)$ satisfy the equation
\begin{equation}\label{16}
T(x+y)+S(x-y)=A(x)+B(y), \quad x, y\in X,
\end{equation}
where $A(x)=T(x)+S(x)$, $B(y)=T(y)+S(-y)$. We solve the equation  (\ref{16}) using the finite difference method. We replace $x$ with $x+k$ and $y$ with $y+k$ in (\ref{16}), where $k$ is an arbitrary element of the group   $X$. Subtract   (\ref{16}) from the obtained equation. We have
\begin{equation}\label{17}
\Delta_{2k}T(x+y)=\Delta_{k}A(x)+\Delta_{k}B(y).
\end{equation}
We replace $y$ with $y+h$ in (\ref{17}), where     $h$ is an arbitrary element of the group   $X$. Subtract (\ref{17}) from the obtained equation. We get
\begin{equation}\label{18}
\Delta_{h}\Delta_{2k}T(x+y)=\Delta_{h}\Delta_{k}B(y).
\end{equation}
We replace $x$ with $x+h$ in (\ref{18}) and   subtract (\ref{18}) from the obtained equation. We have
$$
\Delta_{h}^2\Delta_{2k}T(x+y)=0.
$$
Putting here $y=0$, we receive that the function  $T(x)$  satisfies the equation (\ref{5}). By Lemma \ref{l2}, the function   $T(x)$ is represented in the form  (\ref{6}). This implies that the function $f(x)$ can be represented in the form
$$
 f(x)=\exp\{P(x)+l(x)+r(x)\}, \quad x\in X,
$$
where the function $P(x)$ satisfies the equation    $(\ref{2})$,
  the function $l(x)$ satisfies the equation $(\ref{3})$,
    the function $r(x)$ is  a constant    on each of the cosets of the subgroup $X^{(2)}$ in  $X$.

Reasoning similarly, we obtain that the function    $S(x)$  also satisfies the equation   (\ref{5}) and by Lemma \ref{l2},  the function $g(x)$ may be represented in the form
$$
g(x)=\exp\{Q(x)+m(x)+s(x)\}, \quad x\in X,
$$
where the function $Q(x)$ satisfies the equation    $(\ref{2})$,
  the function $m(x)$ satisfies the equation $(\ref{3})$,
    the function $s(x)$ is  a constant    on each of the cosets of the subgroup $X^{(2)}$ in  $X$.
Substituting the obtained expressions for the functions  $f(x)$ and $g(x)$ into the equation (\ref{1}), we find that $P(x)=Q(x)$, $x\in X$. Moreover, assuming that in the expression obtained after substitution
 $x, y\in x_\alpha+X^{(2)}$, we get that $s(x)=-r(x)$, $x\in X$.   $\blacksquare$
 \begin{corollary}\label{c2}
If we assume in Theorem $\ref{th1}$ that the group $X$ satisfies the condition $X^{(2)}=X$, then the functions $f(x)$ and $g(x)$ may be represented in the form
$$
 f(x)=\exp\{P(x)+l(x)+r\}, \quad  g(x)=\exp\{P(x)+m(x)-r\}, \quad x\in X,
$$
where the function $P(x)$ satisfies the equation    $(\ref{2})$,
  each of the functions $l(x)$ and $m(x)$ satisfies the equation $(\ref{3})$, $r$ is a constant.
\end{corollary}
Obviously, for the group of real numbers   $X=\mathbb{R}$ Theorem A follows from Corollary $\ref{c2}$.
\begin{corollary}\label{c3}
If we assume in Theorem $\ref{th1}$ that $f(x)=g(x)$, i.e. instead of the equation $(\ref{1})$ we consider the equation
\begin{equation}\label{n4}
f(x+y)f(x-y)=f^2(x)f(y)f(-y), \ x, y\in X,
\end{equation}
then the function  $f(x)$   may be represented in the form
$$
 f(x)=\exp\{P(x)+l(x)\}, \quad    x\in X,
$$
 where the function $P(x)$ satisfies the equation    $(\ref{2})$ and
  the function $l(x)$ satisfies the equation $(\ref{3})$.
\end{corollary}
 \begin{remark}\label{r1}
 The description of the solutions of the equation (\ref{1}) given in Theorem $\ref{th1}$ is complete in the following sense.  If  functions $f(x)$ and $g(x)$  are represented in the form (\ref{rev1}), then they satisfy the equation
 $(\ref{1})$. To verify this, it is enough to note that for any $x, y\in X$  the elements $x+y$ and $x-y$ are in the same coset  of the subgroup $X^{(2)}$ in  $X$.  It is obvious that the similar statement also holds for Corollaries    \ref{c2} and    \ref{c3}.
\end{remark}

\section{Solution of the Kac--Bernstein functional equation in the class of complex-valued functions satisfying  the Hermitian condition}

Let $X$ be an Abelian group. We say that a complex-valued function $f(x)$ on the group $X$ satisfies the Hermitian condition if
\begin{equation}\label{19}
f(-x)=\overline{f(x)}, \quad x\in X.
\end{equation}
The characteristic functions of probability distributions on a locally compact Abelian group  satisfy, in particular, the Hermitian condition.
We note that in the article   \cite{TVP2016} the solutions of the equation   (\ref{1}) in the class of characteristic functions on such groups were studied. In this section we study the solutions of the equation   (\ref{1}) in the class of functions satisfying the Hermitian condition on an arbitrary group $X$. We complete  Theorem $\ref{th1}$ with the following statement.
\begin{theorem}\label{th2}
Let $X$ be an Abelian group. Let $f(x)$ and $g(x)$  be complex-valued non-vanishing functions on the group  $X$, each of which satisfies  the Hermitian condition $(\ref{19})$.
If the functions $f(x)$ and $g(x)$  satisfy  the Kac--Bernstein functional equation   $(\ref{1})$, then the functions $f(x)$ and $g(x)$ may be represented in the form
\begin{equation}\label{n9}
f(x)=\alpha(x)a(x)\exp\{P(x)+r(x)\}, \quad g(x)=\beta(x)b(x)\exp\{P(x)-r(x)\}, \quad x\in X,
\end{equation}
where each of the complex-valued functions $\alpha(x)$  and $\beta(x)$ satisfies the equation
\begin{equation}\label{20}
\alpha(x+y)=\alpha(x)\alpha(y), \quad x, y\in X,
\end{equation}
and the condition
\begin{equation}\label{21}
|\alpha(x)|=1, \quad x\in X,
\end{equation}
the functions $a(x)$ and $b(x)$ take the values $\pm 1$ and they
are constant on each of the cosets of the subgroup $X^{(4)}$ in  $X$,
the real-valued function $P(x)$ satisfies the equation    $(\ref{2})$,
      the real-valued function $r(x)$  is  a constant    on each of the cosets of the subgroup $X^{(2)}$ in  $X$.
\end{theorem}
{\bf Proof}.  It is clear that without loss of generality   we can assume that $f(0)=g(0)=1$. Since the functions $|f(x)|$ and $|g(x)|$, obviously also satisfy the equation $(\ref{1})$ and they are even, by Theorem \ref{th1},   there is a representation
\begin{equation}\label{22}
|f(x)|=\exp\{P(x)+r(x)\}, \quad |g(x)|=\exp\{P(x)-r(x)\}, \quad x\in X,
\end{equation}
where the   function $P(x)$  is real-valued and satisfies the equation    $(\ref{2})$,
      the function $r(x)$  is real-valued and is a constant    on each of the cosets of the subgroup $X^{(2)}$ in  $X$. Put
    \begin{equation}\label{n7}
p(x)=f(x)/|f(x)|, \quad q(x)=g(x)/|g(x)|, \quad x\in X.
\end{equation}
Then $|p(x)|=|q(x)|=1$ for $x\in X$, and the functions $p(x)$ and $q(x)$ satisfy the equation  $(\ref{1})$ which takes the form
  \begin{equation}\label{23}
p(x+y)q(x-y)=p(x)p(y)q(x)q(-y), \quad x, y\in X.
\end{equation}
Putting $y=x$, then $y=-x$ in (\ref{23}) we obtain
 \begin{equation}\label{24}
p(2x)=p^2(x), \quad  q(2x)=q^2(x), \quad x \in X.
\end{equation}
Interchange $x$ and $y$ in (\ref{23}), and multiply  the resulting equation and the equation (\ref{23}). Taking into account (\ref{19}), we find
\begin{equation}\label{25}
p^2(x+y)=p^2(x)p^2(y), \quad   x, y \in X.
\end{equation}
It follows from (\ref{24}) and (\ref{25}) that the restriction of the function $p(x)$ to the subgroup $X^{(2)}$ satisfies the equation  (\ref{20}). Reasoning similarly we get that the restriction of the function $q(x)$ to the subgroup $X^{(2)}$ also satisfies the equation  (\ref{20}).

The group $X$ can be considered as a discrete Abelian group. It is well known that then every function on an arbitrary subgroup of the group   $X$   satisfying the equation (\ref{20}) and the condition  (\ref{21}) is a restriction to this subgroup a function on the group    $X$,  satisfying on $X$ the equation (\ref{20}) and the condition  (\ref{21}). We apply this result to the subgroup  $X^{(2)}$. Let $\alpha(x)$ and $\beta(x)$ be the corresponding extensions from the subgroup $X^{(2)}$ to the group $X$ of the functions   $p(x)$ and $q(x)$.    Note that the functions $\alpha(x)$ and $\beta(x)$ satisfy the equation (\ref{1}) on $X$.  Put
 \begin{equation}\label{n8}
a(x)=p(x)/\alpha(x), \quad b(x)=q(x)/\beta(x), \quad x\in X.
\end{equation}
Since $a(x)=b(x)=1$ for $x\in X^{(2)}$, it follows from (\ref{24}) that the functions $a(x)$ and $b(x)$ are even and take the values   $\pm 1$.  It is obvious that the functions  $a(x)$ and $b(x)$ satisfy the equation (\ref{23})  which takes the form
\begin{equation}\label{26}
a(x+y)b(x-y)=a(x)a(y)b(x)b(y), \quad x, y\in X.
\end{equation}
If $x, y\in x_\alpha+X^{(2)}$, then $x\pm y\in X^{(2)}$, and the left-hand side of the equation   (\ref{26}) is equal to 1. This implies that  $a(x)b(y)=a(y)b(x)$ for $x, y\in  x_\alpha+X^{(2)}$, and hence either    $a(x)=b(x)$ or $a(x)=-b(x)$ for $x\in  x_\alpha+X^{(2)}$.  Given this, it follows from (\ref{26}) that
\begin{equation}\label{27}
a(x+y)a(x-y)=1, \quad b(x+y)b(x-y)=1, \quad x, y\in  X^{(2)}{\cup}\big(x_\alpha+X^{(2)}\big),
\end{equation}
and (\ref{27}) implies that the functions $a(x)$ and $b(x)$ are    constant   on each of the cosets of the subgroup $X^{(4)}$ in  $X$.  The statement of the theorem follows from    (\ref{22}), (\ref{n7}) and (\ref{n8}).  $\blacksquare$\begin{remark}\label{r3}
The functions $a(x)$ and $b(x)$ arising in the proof of Theorem $\ref{th2}$  are constant on each of the cosets of the subgroup $X^{(4)}$ in $X$. Surprisingly   is the fact that in this statement we cannot replace $X^{(4)}$ by $X^{(2)}$. Below we construct an   example that shows it (compare with Remark 2 in \cite{TVP2016}).

 Denote by $\mathbb{Z}(4)=\{0, 1, 2, 3\}$ the group of  residue  classes modulo 4. Let $X=(\mathbb{Z}(4))^2$.  Denote by $x=(m, n)$, $m, n\in \mathbb{Z}(4)$, elements of the group $X$. Then $X^{(2)}=\{(0, 0), (0, 2), (2, 0), (2, 2)\}$, $X^{(4)}=\{(0, 0)\}.$
Consider on the group $X$ the functions
$$f(m,n)=\begin{cases} 1, \ \ \ (m,n) \in \{X^{(2)}, (1,0), (3,0), (0,1),
(0,3),   (1,1),   (3,3)\}, \\ -1, \   (m,n) \in \{(1,2),   (3,2),   (2,1), (2,3),   (1,3),   (3,1 )\}\\
\end{cases}
$$
and
$$
g(m,n)=\begin{cases} 1, \ \ \ (m,n) \in \{X^{(2)},   (1,0), (3,0), (0,1), (0,3), (1,3),   (3,1)\}, \\ -1, \
 (m,n) \in \{(1,2), (3,2),   (2,1), (2,3),  (1,1),   (3,3 )\}.
\end{cases}
$$
Verify that the functions $f(x)$ and $g(x)$ satisfy the equation (\ref{1}).   Let
$$
X=X^{(2)}\cup\big((1, 0)+X^{(2)}\big)\cup\big((0, 1)+X^{(2)}\big)\cup\big((1, 1)+X^{(2)}\big)
$$
be a decomposition of the group $X$ with respect to the subgroup $X^{(2)}$. It follows from the expressions for the functions $f(x)$ and $g(x)$ that $f(x)g(x)=1$ for $x\in X^{(2)}\cup\big((0, 1)+X^{(2)}\big)\cup\big((1, 0)+X^{(2)}\big)$ and $f(x)g(x)=-1$ for $x\in (1, 1)+X^{(2)}$.

We consider successively the cases when $x$ and $y$ belong to the given cosets s of the subgroup $X^{(2)}$ in  $X$.

Let $x, y\in X^{(2)}$. Then $x\pm y\in X^{(2)}$, and   the right-hand and left-hand sides of the equation (\ref{1}) are equal to $1$.

Let $x \in X^{(2)}$ and $y\in (1, 0)+X^{(2)}$.  Then the right-hand  side of the equation   (\ref{1})  is equal to $1$. We note that $x\pm y\in (1, 0)+X^{(2)}$. If the left-hand side of the equation  (\ref{1}) is equal to $-1$, then either $x+y\in \{(1, 0), (3, 0)\}$,  $x-y\in \{(1, 2), (3, 2)\}$    or $x+y\in \{(1, 2), (3, 2)\}$, $x-y\in \{(1, 0), (3, 0)\}$. In each of these cases  $2x\ne (0, 0)$, which contradicts the fact that   $x \in X^{(2)}$. Hence, the left-hand side of the equation (\ref{1}) is also equal to $1$.

Let $x, y\in (1, 0)+X^{(2)}$. Then the right-hand  side of the equation (\ref{1})  is equal to $1$. We note that $x\pm y\in X^{(2)}$. Hence, the left-hand side of the equation (\ref{1}) is also equal to $1$.

Let $x \in (1, 0)+X^{(2)}$ and $y\in (0, 1)+X^{(2)}$. Then the right-hand  side of the equation (\ref{1})  is equal to $1$.  We note that $x\pm y\in (1, 1)+X^{(2)}$.   If the left-hand side of the equation  (\ref{1}) is equal to $-1$, then either $x\pm y\in \{(1, 1), (3, 3)\}$ or $x\pm y\in \{(1, 3), (3, 1)\}$.    In each of these cases  $2x\in \{(0, 0), (2, 2)\}$, which contradicts the fact that $x \in (1, 0)+X^{(2)}$. Hence, the left-hand side of the equation (\ref{1}) is also equal to $1$.

Let   $x\in (1, 0)+X^{(2)}$ and $y\in (1, 1)+X^{(2)}$. Then the right-hand  side of the equation (\ref{1})  is equal to $-1$.  We note that $x\pm y\in  (0, 1)+X^{(2)}$. If the left-hand side of the equation  (\ref{1}) is equal to $1$, then either $x\pm y\in \{(0, 1), (0, 3)\}$ or $x\pm y\in \{(2, 1), (2, 3)\}$. In each of these cases  $2x\in \{(0, 0), (0, 2)\}$, which contradicts the fact that $x\in (1, 0)+X^{(2)}$. Hence, the left-hand side of the equation (\ref{1}) is also equal to $-1$.

The remaining cases are similar to those considered above.

In the constructed example in the representation  (\ref{n9}) $\alpha(x)\equiv 1$, $\beta(x)\equiv 1$, $P(x)\equiv 0$ and $r(x)\equiv 0$ for $x\in X$,  but the functions $a(x)=f(x)$ and $b(x)=g(x)$ are not constant on each of the cosets   of the subgroup $X^{(2)}$ in  $X$.
\end{remark}
\begin{remark}\label{r5}
The example in Remark \ref{r3} also shows that both  the function $a(x)$ and the function  $b(x)$ arising in the proof of Theorem $\ref{th2}$ need not   satisfy the equation (\ref{20}). At the same time, assuming that the group $X$ is such that   $X/X^{(2)}\cong \mathbb{Z}(2)$, where $\mathbb{Z}(2)$ is the group of  residue  classes modulo 2, each of the functions $a(x)$ and $b(x)$ satisfies the equation  (\ref{20}). Indeed,
if $X/X^{(2)}\cong \mathbb{Z}(2)$, then a decomposition of the group $X$ with respect to the subgroup $X^{(2)}$ is of the form   $X=X^{(2)}\cup(x_1+X^{(2)})$. Therefore, (\ref{27}) implies that the functions   $a(x)$ and  $b(x)$  are     constant    on each of the cosets of the subgroup $X^{(2)}$ in  $X$. Since $a(x)=b(x)=1$ for $x\in X^{(2)}$ and both functions    $a(x)$  and   $b(x)$  on the coset $x_1+X^{(2)}$ take the values $\pm 1$, then, obviously,   each of them satisfies the equation (\ref{20}). The representation (\ref{n9}) then takes the form
$$
f(x)=\pm\alpha(x)\exp\{P(x)+r(x)\}, \quad g(x)=\pm\beta(x)\exp\{P(x)-r(x)\}, \quad x\in X.
$$

If we assume that the group  $X$ satisfies the condition $X^{(2)}=X$, then    Theorem \ref{th2} implies directly that in this case the representation (\ref{n9})   takes the form
\begin{equation}\label{02.01.2}
f(x)=\pm\alpha(x)\exp\{P(x)+r\}, \quad g(x)=\pm\beta(x)\exp\{P(x)-r\}, \quad x\in X,
\end{equation}
where $r$ is a real constant.
\end{remark}
 \begin{remark}\label{r4} Generally speaking,  if  functions $f(x)$  and $g(x)$ are represented in the form (\ref{n9}), then they need not satisfy the equation (\ref{1}). An example can be easily constructed on the group $(\mathbb{Z}(4))^2$. Nevertheless, if we additionally require that the functions   $a(x)$ and $b(x)$  in the representation  (\ref{n9}) are   constant on each of the cosets   of the subgroup $X^{(2)}$ in  $X$ and $a(x)b(x)=1$ for $x\in X$, it is easy to verify that in this case the functions $f(x)$  and $g(x)$   already satisfy the equation (\ref{1}).
\end{remark}
\begin{corollary}\label{c5}
If we assume in Theorem $\ref{th2}$ that $f(x)=g(x)$, i.e. instead of the equation $(\ref{1})$ we consider the equation   $(\ref{n4})$,
then the function  $f(x)$   may be represented in the form
 \begin{equation}\label{n5}
f(x)=\alpha(x)a(x)\exp\{P(x)\}, \quad x\in X,
\end{equation}
where the complex-valued function  $\alpha(x)$    satisfies the equation $(\ref{20})$ and the condition  $(\ref{21})$, the function  $a(x)$  takes the values $\pm 1$ and is  a  constant   on each of the cosets of the subgroup $X^{(2)}$ in  $X$ and the real-valued function $P(x)$ satisfies the equation $(\ref{2})$.
\end{corollary}
{\bf Proof}. Following the scheme of the proof of Theorem \ref{th2}, we arrive at the equation   $(\ref{26})$. Since $f(x)=g(x)$, we have $a(x)=b(x)$ in (\ref{26}) and we get that the function   $a(x)$ satisfies the equation
\begin{equation}\label{02.01.1}
a(x+y)a(x-y)=1, \quad x, y\in X.
\end{equation}
It follows from (\ref{02.01.1}) that the function $a(x)$      is  a  constant   on each of the cosets of the subgroup $X^{(2)}$ in  $X$. $\blacksquare$\begin{remark}\label{r2}
Despite the fact that the function $a(x)$ in the representation  (\ref{n5}) in Corollary \ref{c5} is  a  constant   on each of the cosets of the subgroup $X^{(2)}$ in  $X$, unlike the representation  (\ref{n9}) in Theorem \ref{th2},  where  the function $a(x)$ is  a  constant   on each of the cosets of the subgroup $X^{(4)}$ in  $X$, generally speaking the function $a(x)$ need  not satisfy the equation $(\ref{20})$.
We construct an example. Denote by $\mathbb{Z}$ the group of integers. Let  $X=\mathbb{Z}^2$. Denote by $x=(m, n)$, $m, n\in \mathbb{Z}$, elements of the group $X$. Put
$$
f(m, n)=\exp\{i\pi mn\}, \quad (m, n)\in X.
$$
It is obvious that the function  $f(x)$ satisfies the equation (\ref{n4}). It is easy to see that in this case  in the representation (\ref{n5}) $\alpha(x)\equiv 1$, $P(x)\equiv 0$  for $x\in X$,  but the function $a(x)=f(x)$ does not satisfy the equation $(\ref{20})$.
\end{remark}
\begin{remark}\label{r01}
In contrast to Theorem $\ref{th2}$, the description of the solutions of the equation (\ref{n4}) given in Corollary  $\ref{c5}$ is complete in the following sense. If  a function  $f(x)$ is represented in the form (\ref{n5}), then it satisfies the equation (\ref{n4}). It follows from the fact that if a function  $a(x)$  takes the values $\pm 1$ and is  a  constant   on each of the cosets of the subgroup $X^{(2)}$ in  $X$, then it satisfies the equation   (\ref{02.01.1}), and hence also satisfies the equation (\ref{n4}).
\end{remark}
Assume that an Abelian group $X$ satisfies the condition $X^{(2)}=X$. Based on the representation (\ref{02.01.2}), we prove the following statement.
\begin{theorem}\label{th3}
Let $X$ be an Abelian group, such that $X^{(2)}=X$.  Let   $f(x)$ and $g(x)$  be complex-valued functions that are not identically zero on the group $X$, each of which satisfies  the Hermitian condition $(\ref{19})$.
If the functions $f(x)$ and $g(x)$  satisfy  the Kac--Bernstein functional equation   $(\ref{1})$, then they   can be represented in the form
\begin{equation}\label{rev2}
f(x)=\begin{cases} \pm\alpha(x)\exp\{P(x)+r\}, \    x\in G,
\\ 0, \qquad\qquad\qquad\quad\qquad    x \notin G,
\end{cases} \quad
g(x)=\begin{cases} \pm\beta(x)\exp\{P(x)-r\}, \    x\in G,
\\ 0, \qquad\qquad\qquad\quad\qquad    x \notin G,
\end{cases}
\end{equation}
where $G$ is a subgroup of $X$, such that the factor-group  $X/G$ contains no elements of order   $2$,  each of the complex-valued functions $\alpha(x)$  and $\beta(x)$ satisfies the equation $(\ref{20})$
and the condition $(\ref{21})$,
 the real-valued function $P(x)$ satisfies the equation   $(\ref{2})$, $r$ is a real constant.
\end{theorem}
{\bf Proof}. It follows from the equation (\ref{1}) that either  $f(0)g(0)=0$ or   $f(0)g(0)=1$. Let  $f(0)g(0)=0$ and assume for definiteness that  $f(0)=0$. Then, putting $y=-x$ in the equation (\ref{1}) we obtain that $|f(x)|^2g^2(x)=0$ for $x\in X$. Hence, $f(x+y)g(x-y)=0$ for $x, y\in X$. In view of   $X^{(2)}=X$, this implies that either   $f(x)\equiv 0$ or $g(x)\equiv 0$ for $x\in X$, contrary to the condition. Therefore $f(0)g(0)=1$. Without loss of generality, we can assume that $f(0)=g(0)=1$.

Putting $y=x$, then $y=-x$ in (\ref{1}) we obtain
\begin{equation}\label{28}
f(2x)=f^2(x)|g^2(x)|, \quad   g(2x)=|f^2(x)|g^2(x), \quad x\in X.
\end{equation}
In view of $X^{(2)}=X$, it follows from (\ref{28})  that
\begin{equation}\label{29}
|f(x)|=|g(x|, \quad   x\in X.
\end{equation}
Taking into account (\ref{29}), the equation (\ref{1}) implies that  the set
$$
G=\{x\in X: f(x)\ne 0\}=\{x\in X: g(x)\ne 0\}
$$
is a subgroup of $X$. It follows from (\ref{28}) that the subgroup  $G$ has the property: if $2x\in G$, then $x\in G$. This implies that the factor-group $X/G$ contains no elements of order $2$. It also follows from  (\ref{28}) that the subgroup $G$ satisfies the condition  $G^{(2)}=G$. We apply the representation   (\ref{02.01.2})  to the group $G$.  $\blacksquare$

\begin{remark}\label{r6}
The description of the solutions of the equation (\ref{1}) given in Theorem  $\ref{th3}$ is complete in the following sense. Let $X$ be an Abelian group, such that $X^{(2)}=X$.  If   functions  $f(x)$ and $g(x)$ are represented  in the form (\ref{rev2}), then they satisfy the equation (\ref{1}).
\end{remark}

 \end{document}